\documentclass{article}
\usepackage[utf8]{inputenc}

\usepackage[utf8]{inputenc}
\usepackage[T1]{fontenc}
\usepackage{lmodern}

\usepackage{amsmath}
\usepackage{mathtools}
\usepackage{amssymb}
\usepackage{amsthm}

\usepackage{todonotes}

\usepackage{array}

\makeatletter
\def\@fnsymbol#1{\ensuremath{\ifcase#1\or \star\or \ddagger\or
   \mathsection\or \mathparagraph\or \|\or **\or \dagger\dagger
   \or \ddagger\ddagger \else\@ctrerr\fi}}
\makeatother

\usepackage[plainpages=false,pdfpagelabels,colorlinks=true,citecolor=blue,hypertexnames=false]{hyperref}

\newtheorem{theorem}{Theorem}
\newtheorem{proposition}[theorem]{Proposition}

\theoremstyle{definition}
\newtheorem{remark}[theorem]{Remark}

\newcommand{\ps}[1]{[\![#1]\!]}

\newcommand{\CC}{\mathbb{C}}
\newcommand{\NN}{\mathbb{Z}}
\newcommand{\ZZ}{\mathbb{Z}}
\newcommand{\qbinom}{\genfrac{[}{]}{0pt}{}}
\DeclarePairedDelimiter\floor{\lfloor}{\rfloor}

\title{On the $q$-analogue of Pólya's Theorem}
\author{Alin Bostan\footnote{Inria, Univ. Paris-Saclay, France, \url{alin.bostan@inria.fr}.} 
\ and Sergey Yurkevich\footnote{U. Wien, Austria and Inria, Univ. Paris-Saclay, France, \url{sergey.yurkevich@univie.ac.at}.}}

\begin{document}

\maketitle

\begin{abstract}
    We answer a question posed by Michael Aissen in 1979 
    about the $q$-analogue of a classical theorem of George Pólya (1922) on the algebraicity of
    (generalized) diagonals 
    of bivariate rational power series. In particular, we prove that the answer to Aissen's question, in which he considers $q$ as a variable, is negative in general. Moreover, we show that when $q$ is a complex number, the answer is positive if and only if $q$ is a root of unity.
\end{abstract}

\section{Introduction}\label{sec:1}

A beautiful but rather unknown theorem of Pólya~\cite{Polya22} states the following: 
\begin{quote}
Given two algebraic power series\footnote{Recall that $f(x) \in \CC\ps{x}$ is called \emph{algebraic} if there exists a bivariate non-zero polynomial $P(x,z)$ in $\CC[x,z]$ such that $P(x,f(x)) = 0$. A non-algebraic series is called \emph{transcendental}.} $\varphi(x)$ and $\Phi(x)$, let $A_{(i,j)}$ be the coefficient of $x^j$ in $\Phi(x)\varphi(x)^i$. Consider a straight line in the plane and let $(p_n)_{n\geq0}$ be the sequence of non-negative integer lattice points in~$\ZZ^2$ lying on this line. Then $F(x) = \sum_{n \geq 0} A_{p_n}x^n$ is algebraic. 
\end{quote}

In particular, this theorem implies that the generalized diagonal $\Delta_{a,b}$ of a bivariate rational power series is algebraic, where for $f(x,y) = \sum_{i,j \geq 0 } f_{i,j}x^iy^j$, we define $\Delta_{a,b}(f) \coloneqq \sum_{n \geq 0} f_{a n,b n} x^n$. For example, one finds
\begin{align*}
    \Delta_{1,1}\left( \frac{1}{1-x-y}\right) = \Delta_{1,1}\left( \sum_{i,j \geq 0 } \binom{i+j}{i} x^i y^j \right) = \sum_{n \geq 0} \binom{2 n}{n} x^n =  \frac{1}{\sqrt{1-4x}},
\end{align*}
and the latter is a root of $P(x,z) = (1-4x)z^2-1$.

In Pólya's formulation, this example is obtained by choosing $\Phi(x)=1$, $\varphi(x) = 1+x$ and the main diagonal $\{ x=y \}$ of $\mathbb{Z}^2$. In fact, this special case is the main foundation for the observation and question which led to this article.

There is a more combinatorial rephrasing of this example. Arrange Pascal's triangle in the following way:
\begin{center}
\begin{tabular}{>{$}l<{$}|*{7}{c}}
0 &1&&&&&&\\
1 &1&1&&&&&\\
2 &1&2&1&&&&\\
3 &1&3&3&1&&&\\
4 &1&4&6&4&1&&\\
5 &1&5&10&10&5&1&\\
6 &1&6&15&20&15&6&1\\
\vdots& &&\dots&\dots&&&\\
\hline
\multicolumn{1}{l}{} &0&1&2&3&4&5&6\\
\end{tabular}
\end{center}
\noindent
and consider a line passing through infinitely many lattice points of the above triangle. If we denote the resulting sequence of values on these lattice points by $(u_j)_{j\geq0}$, then Pólya's theorem ensures that the generating function $f(x)=\sum_{j\geq0} u_{j}x^j$
is algebraic.

It is easy to see that the above condition on the line can be reformulated into the existence of non-negative integers $n,k,a,b$ with $n \geq k, a \geq b$, $\gcd(a,b) =1$ and that either $n-a<k-b$ or $k-b<0$, such that $u_j = \binom{n+aj}{k+bj}$. So this special case of Pólya's theorem simply asserts that
\[
\sum_{j\geq 0} \binom{n+aj}{k+bj}x^j \in \CC\ps{x} \quad \text{is algebraic.}
\]

It is well-known and easy to see that the binomial coefficient $\binom{x+y}{x}$ counts lattice paths from the origin to $(x,y) \in \ZZ^2$ with only North and East steps. The observation above therefore implies that the generating function of the number of such paths, as $(x,y)$ increases on a line, is algebraic. Now recall the definition of the $q$-analogue of the binomial coefficient, called the \emph{$q$-binomial} coefficient:
\begin{align*}
    \qbinom{n}{k}_q &\coloneqq \frac{[n]_q!}{[k]_q![n-k]_q!}, \qquad \text{where}\\
    [n]_q! &\coloneqq (1+q) \cdots (1+q+\cdots +q^{n-1}),
\end{align*}
for a variable $q$ and integers $n,k$ with $0 \leq k\leq n$. It is not difficult to check that $\qbinom{n}{k}_q \in \ZZ[q]$ is a polynomial in $q$ of degree $k(n-k)$. 
Arranged as in the figure above, these $q$-binomial coefficients give rise to the so-called $q$-Pascal triangle.

Pólya showed \cite{Polya69} that the coefficient of $q^j$ in $\qbinom{x+y}{x}_q$ counts lattice paths in~$\ZZ^2$ from the origin to $(x,y)$ with same steps as before and with area underneath equal to $j$, see also \cite{Aissen79}.
Aissen asked in \cite[p.~585]{CaMaMcMo79} the natural question whether the following $q$-analogue of Pólya's statement about algebraicity of such path generating functions holds:
\begin{quote}
    Fix integers $n,k,a,b$ with $n \geq k \geq 0, a \geq b \geq 0$ and $\gcd(a,b) =1$. Moreover assume that either $n-a<k-b$ or $k-b<0$. Let 
    \[
        F(x,q) \coloneqq \sum_{j\geq 0} \qbinom{n+aj}{k+bj}_qx^j
        \in \CC[q]\ps{x}.
    \]
    Is the power series $F(x,q)$ algebraic? That is, does there exist a non-zero polynomial $P(x,q,z) \in \CC[x,q,z]$ such that $P(x,q,F(x,q)) = 0$?
\end{quote}

The inequality conditions $n \geq k \geq 0$ and $a \geq b \geq 0$ ensure that $F(x,q)$ is well-defined and not a polynomial. The condition $\gcd(a,b)=1$ means that the line passing through the $q$-Pascal triangle does not ``skip'' terms; as we will see, it does not affect the algebraicity of $F(x,q)$. Also the condition $n-a<k-b$ or $k-b<0$ is just a translation of the geometric picture that $F(x,q)$ collects all terms on that line. Aissen noticed the following fact: if the line is parallel to an edge of the $q$-Pascal triangle (i.e., if $a = b$ or $b = 0$), then $F(x, q)$ is trivially algebraic, because it is actually a rational function of $x$ and $q$. Hence, in what follows, we will assume that $a \neq b$ and $b \neq 0$. More precisely, we will only consider 
\emph{admissible integers} $n,k,a,b$, in the following sense: $n \geq k \geq 0$, $a > b > 0$, $\gcd(a,b) =1$, and either 
$n-k<a-b$ or $k<b$.

\section{Results}

Using elementary asymptotic estimates on the coefficients, we can show that the answer to Aissen's question is negative, because for any fixed $z \in \CC$ of absolute value larger than $1$ the coefficients of $F(x,z)$ grow too fast for this series to be algebraic. However, we notice that the same argument does not apply for the univariate power series $h_z(x) \coloneqq F(x,z)$, where $z \in \CC$ with $|z| \leq 1$. In fact, our main result (Theorem~\ref{thm:main} below) is that $h_z(x)$ is algebraic if and only if $z$ is a root of unity. 

We start with a particular case of our main result and the answer to Aissen's question: the generating function of the central $q$-binomial coefficients is not algebraic.

\begin{proposition}\label{F_not_algebraic}
For $n=k=0$ and $a=2,b=1$ the series $F(x,q)$ is not algebraic.
\end{proposition}
\begin{proof}
Assume that $F(x,q)$ is algebraic with minimal polynomial $P(x,q,z)$. Then the series $h_2(x) \coloneqq F(x,2)$ must be algebraic as well, since $P_2(x,z) \coloneqq P(x,2,z) \not \equiv 0 $ satisfies $P_2(x,h_2(x)) = 0$. We have 
\begin{align*}
    h_2(x) &= \sum_{j\geq 0}\frac{(2^{j+1}-1)(2^{{j+2}}-1)\cdots (2^{{2j}}-1)}{(2-1)(2^{2}-1)\cdots (2^{j}-1)}x^j \\
    & = 1 + 3x + 35x^2 + 1395x^3 + 200787x^4 + \cdots.
\end{align*}
Using the obvious inequality $(2^{j+k}-1)/(2^k-1) > 2^j$, we see that the $j$-th coefficient of $h_2(x)$ is greater than $2^{j^2}$. This growth rate is too fast for $h_2(x)$ to be algebraic, see e.g., Theorem D in~\cite{Flajolet87}, or Theorem 3 in~\cite{RiRo14}.
\end{proof}

A more elementary way to see that the growth rate of the coefficients of $h_2(x)$ is incompatible with algebraicity of the function is to notice that this rate is too fast even for D-finite functions. Recall~\cite{Stanley80} that a power series $f(x) = \sum_{i\geq0} u_ix^i$ is called \emph{D-finite}  if it satisfies a linear differential equation with polynomial coefficients:
\[
p_n(x) f^{(n)}(x) + \cdots + p_0(x) f(x) = 0. 
\]
A classical theorem ensures that any algebraic function is D-finite~\cite[Thm.~2.1]{Stanley80} (see also~\cite{BCLSS07}), whereas the latter class of functions is clearly much larger. An equivalent characterization of D-finite  series~\cite[Thm.~1.5]{Stanley80} states that the coefficients sequence satisfies a linear recurrence with polynomial coefficients:
\[
 u_{j+r}c_r(j) + \dots + u_{j}c_0(j) = 0 , \quad j \geq 0.
\]
A sequence $(u_j)_{j\geq 0}$ is called \emph{P-recursive}  if it satisfies a recurrence as above. 

A simple estimation on the growth rate of such sequences shows that any P-recursive sequence $(u_j)_{j\geq 0}$ grows at most like a power of $j!$ which is slower than $2^{j^2}$. We will use the fact that the coefficient sequence of an algebraic function is necessarily P-recursive again later. \\

We have just proved that the bivariate series $F(x,q)$ cannot be algebraic for all admissible $n,k,a,b$. In the same manner, we can prove that $F(x,q)$ is not algebraic for \emph{any} admissible $n,k,a,b$. This is easily reduced to the following task:
when is the univariate power series $h_z(x) \coloneqq F(x,z)$ (for a fixed $z\in\CC\setminus \{ 0 \}$) algebraic? Obviously, the same argument as in the proof of Proposition~\ref{F_not_algebraic} applies for any $z$ with $|z|>1$: the growth rate of $z^{j^2}$ is incompatible with algebraicity. However, for $z=\omega$ on the unit circle or for $|z|<1$ the same argument does not work; in fact we have the following result:

\begin{proposition} \label{h_algebraic}
    Let $n,k,a,b$ be admissible integers and $\omega \in \CC$ be a root of unity. Then $h_\omega(x)$ is algebraic.
\end{proposition}

Before proving it, we recall the $q$-Lucas theorem (see~\cite[p.~22]{Desarmenien82} for an algebraic proof and \cite[p.~131--132]{Sagan92} for a combinatorial proof). It is the $q$-analogue of the well-known Lucas theorem, one form of which states that
\[
\binom{n}{m} \equiv \binom{\floor{\frac{n}{p}}}{\floor{\frac{m}{p}}}\cdot \binom{n-p\floor{\frac{n}{p}}}{m-p\floor{\frac{m}{p}}} \mod p,
\]
where $p$ is a prime number and $m,n$ are non-negative integers; here $\floor{x}$ denotes the integer part of~$x$, that is the largest integer at most equal to~$x$.

\begin{theorem}[$q$-Lucas Theorem]
    Let $x,y$ be non-negative integers and $\omega \in \CC$ be a root of unity of order $s$. Then
    \[
        \qbinom{x}{y}_\omega = \binom{\floor{\frac{x}{s}}}{\floor{\frac{y}{s}}}\cdot\qbinom{x-s\floor{\frac{x}{s}}}{y-s\floor{\frac{y}{s}}}_\omega.
    \]
\end{theorem}
\noindent
Now we can show that for roots of unity $\omega \in \CC$, the series $h_\omega(x)$ is algebraic.
\begin{proof}[Proof of Proposition \ref{h_algebraic}]
    Let $s$ be the order of $\omega$. We have
    \begin{align*}
        h_\omega(x) = \sum_{j\geq 0} \qbinom{n+aj}{k+bj}_\omega x^j = \sum_{r=0}^{s-1} \sum_{\substack{j \geq 0 \\ j \equiv r \textrm{ mod } s}} \qbinom{n+aj}{k+bj}_\omega x^j.
    \end{align*}
    Let us examine the $s$ summands separately using the $q$-Lucas theorem: 
    \begin{align*}
        \sum_{\substack{j \geq 0 \\ j \equiv r \textrm{ mod } s}} \qbinom{n+aj}{k+bj}_\omega x^j & = \sum_{\ell \geq 0} \qbinom{n+a(\ell s+r)}{k+b(\ell s+r)}_\omega x^j \\
        & = \sum_{\ell \geq 0} \binom{\floor{\frac{n+ar}{s}} + a\ell}{\floor{\frac{k+br}{s}} + b\ell}\qbinom{n+ar - s\floor{\frac{n+ar}{s}}}{k+br - s\floor{\frac{k+br}{s}}}_\omega x^j \\
        & = \qbinom{n+ar - s\floor{\frac{n+ar}{s}}}{k+br - s\floor{\frac{k+br}{s}}}_\omega \cdot \sum_{\ell \geq 0} \binom{\floor{\frac{n+ar}{s}} + a\ell}{\floor{\frac{k+br}{s}} + b\ell} x^j.
    \end{align*}
    By Pólya's theorem (see Section \ref{sec:1}) the last sum is an algebraic series, and therefore obviously the whole last expression is also algebraic. Then, $h_\omega(x)$ is the sum of $s$ algebraic power series, hence it is algebraic as well.
\end{proof}

\noindent
In the remaining part of the article we will show that if $0 < |z|\leq 1$ but $z$ is not a root of unity, then $h_z(x)$ cannot algebraic. This will prove our main theorem:
\begin{theorem}\label{thm:main}
    Let $n,k,a,b$ be admissible integers and let $q\in\CC\setminus \{ 0 \}$.
    Then $h_q(x)$ is an algebraic power series if and only if $q$ is a root of unity.
\end{theorem}

A natural approach to prove this theorem is to use a result originating from Ramis' work~\cite{Ramis92}, see also~\cite[Corollary~2]{ScSi19}, or~\cite[Theorem~7.1]{BeBo92}. It says that if a function $f$ is algebraic and at the same time satisfies a linear $q$-difference equation, that is, an equation of the form
\[
    f(q^mx) + b_{m-1}(x) f(q^{m-1}x) + \cdots + b_0(x) f(x) = 0
\]
for some rational functions $b_0, \ldots, b_m \in \CC(x)$ not all zero and $q \in \CC\setminus \{ 0 \}$ not a root of unity, then $f$ is actually a rational function. Clearly, $h_q(x)$ satisfies a linear $q$-difference equation. Hence, the result described above ensures that if $h_q(x)$ is algebraic, it must already be a rational function. One can prove that for any non-zero $q$, $h_q(x)$ is rational if and only if $a=b$ or $b=0$. Then, altogether, these facts imply Theorem~\ref{thm:main}.

It turns out that a simple modification of the proof that $h_q(x)$ is never rational for admissible integers already implies a much more general fact which does not require the theory of $q$-difference equations in order to prove our main theorem. More precisely, we will prove directly the following stronger result. 
\begin{theorem} \label{not_dfinite}
    Let $n,k,a,b$ be admissible integers and let $q\in\CC\setminus \{ 0 \}$. Then $h_q(x)$ is D-finite if and only if $q$ is a root of unity.
\end{theorem}
\noindent
Recall that an algebraic function is always D-finite, therefore Theorem~\ref{not_dfinite} along with Proposition~\ref{h_algebraic} will allow us to conclude the validity of Theorem~\ref{thm:main}.

We will make use of the following elementary proposition, that 
we suspect to be well-known although we could not locate 
it in the literature.
\begin{proposition}\label{prop:jqj}
    Let $p(x,y) \in \CC[x,y]$ be a bivariate polynomial and assume that for some $q\in\CC\setminus \{0\}$ not a root of unity, we have $p(j,q^j) = 0$ for all $j \in \mathbb{N}$. Then $p(x,y)=0$.
\end{proposition}
\begin{proof}
    We distinguish three cases: $|q|<1$, $|q|=1$ and $|q|>1$. For the first case, write $p(x,y) = p_0(x) + r(x,y)y^n$ for some natural number $n$ and $p_0(x) = p(x,0) \in \CC[x]$, $r(x,y) \in \CC[x,y]$ such that $r(x,0) \neq 0$. It follows that 
    \[
    0 = p(j,q^j) = p_0(j) + r(j,q^j)q^{nj}.
    \]
    Since $|q|<1$, we must have $r(j,q^j)q^{nj} \to 0$ as $j \to \infty$. Therefore, $\lim_{j\to\infty} p_0(j) =0$ and we obtain that $p_0(x) = 0$. Hence, $r(j,q^j) = 0$ for $j\geq0$. But $r(x,y) = r_0(x) + s(x,y)y$ for some polynomial $s(x,y)$ and non-zero $r_0(x)$. By the same argument, $\lim_{j\to\infty}r_0(j)=0$, however this contradicts $r_0(x)\neq0$.
    
    If $|q|=1$, write $p(x,y) = p_0(y) + p_1(y)x+\cdots+p_d(y)x^d$ for some natural number $d$ and polynomials $p_0(y),\dots,p_d(y) \in \CC[y]$, such that $p_{d}(y) \neq 0$. We have 
    \begin{equation} \label{eq:jqj}
        |p_d(q^j)|j^d = \left| \sum_{k=0}^{d-1} p_k(q^j) j^{k} \right|.
    \end{equation}
    Now the idea is that for some sequence $(j_n)_{n \geq 0}$, the terms $|p_0(q^{j_n})|,\dots,|p_{d-1}(q^{j_n})|$ can be bounded by a constant from above and $|p_d(q^{j_n})|$ is bounded from below by a non-zero constant -- this contradicts (\ref{eq:jqj}) because the left-hand side becomes too large. More precisely, choose $\xi$ on the unit circle which is not a root of $p_d(x)$. Then there exists $\varepsilon > 0$ such that $|p_{k}(\xi)| < 1/\varepsilon$ for all $k=0,\dots,d$ and also $\varepsilon < |p_{d}(\xi)| $. Moreover, since $q$ is not a root of unity, Jacobi's Theorem implies that the set $\{q^{j}| j\in\mathbb{N}\}$ is dense on the unit circle~\cite[Thm~3.13]{Devaney89} (see also \cite{Denjoy32}), consequently there are infinitely many $j$ such that $q^j$ is arbitrarily close to $\xi$. Henceforth, there also exist infinitely many $j$ for which $|p_{k}(q^j)| < 1/\varepsilon$ for all $k=0,\dots,d-1$ and $\varepsilon < |p_{d}(q^j)|$. However, at the sequence of these $j$, this contradicts (\ref{eq:jqj}) since then the left-hand side grows at least like $j^{d} \varepsilon$ and the right-hand side is bounded by $dj^{d-1}/\varepsilon$.
    
    Finally, if $|q|>1$, write $p(x,y) = r_0(x) + r_1(x)y + \cdots + r_n(x)y^n$ for polynomials $r_0(x),\dots,r_n(x) \in \CC[x]$ and some natural number $n$. Clearly, if $p(x,y)$ is non-zero, $n$ must be positive. But then we have
    \[
    |r_0(j) + r_1(j)q^j + \cdots + r_{n-1}(j)q^{(n-1)j}| \leq c j^m |q|^{(n-1)j},
    \]
    for some constants $c,m>0$. This contradicts
    \[
    r_n(j)q^{nj} = - r_0(j) - r_1(j)q^j - \cdots - r_{n-1}(j)q^{(n-1)j}
    \]
    for big enough $j$ and finishes the proof. 
\end{proof}

\begin{remark}
    An alternative, purely algebraic, proof of Proposition~\ref{prop:jqj} follows from the fact that for any $d \geq 0$, writing $D = \binom{d+2}{2}$, the determinant of the $D \times D$ matrix 
    \[
        M(z) = (n^i  z^{n j})_{\substack{0 \leq i+j \leq d,\\ 1 \leq n \leq D}}
    \] 
    is given by a constant times a power of $z$ and a product of cyclotomic polynomials in~$z$. More precisely, assuming that the total degree of $p(x,y) = \sum_{i,j} c_{i,j} x^i y^j$ is $d$, the equations $p(j,q^j) = 0$ for $j=1,\dots,D$ yield the following linear system of equations for the vector of unknowns $c_{i,j}$:
    \[
    M(q) \cdot (c_{i,j})_{0\leq i+j \leq d} = 0.
    \]
    To see why $\det M(z)$ only vanishes for $z$ a root of unity, it is useful to write $M(z) = N(1,z,\dots,z^{d})$, where $N(z_0,\dots,z_{d}) = (n^i  z_j^{n})_{\substack{0 \leq i+j \leq d,\\ 1 \leq n \leq D}}$. Then it remains to prove that $\det N(z_0,\dots,z_{d})$ is a constant times a product of $z_i$'s times a product of $(z_i-z_j)$ for $i \neq j$. This follows from the observation that the transpose of $N$ is a generalized Vandermonde matrix; more precisely, it is a matrix corresponding to the linear map from the space of polynomials with no constant term and degree at most $D$ to $\CC^{D}$ given by
    \[
    P(x) \mapsto (P(z_0),\dots,P(z_d),\vartheta P(z_0),\dots,\vartheta P(z_{d-1}),\dots,\vartheta^d P(z_0)),
    \]
    where $\vartheta = x \frac{x}{\mathrm{d}x}$ is the Euler derivative. After a change of basis from the monomials to the the basis $1, (x-z_0), (x-z_0)(x-z_1), \dots, (x-z_0)^{d-1}(x-z_1)^{d-1}(x-z_2)^{d-2}\cdots(x-z_{d})$ the matrix becomes lower-triangular and the determinant evaluation follows trivially. This purely algebraic proof shows that in Proposition~\ref{prop:jqj} one can replace ``for all $j \in \NN$'' by the weaker condition ``for $j = 1,\dots,(d+1)(d+2)/2$, where $d$ is the total degree of $p(x,y)$''. The proof also shows that the conclusion of the proposition holds as well if $q \in \CC\setminus\{0\}$ is assumed not to be a root of unity of order at most $d$.
\end{remark}

Note that Proposition~\ref{prop:jqj} immediately implies that the function $f(x) = q^x$ is transcendental for any non-zero $q \in \CC$ which is also not a root of unity, because an annihilating polynomial $P(x,z)$ would need to satisfy $P(x,q^x) = 0$ and hence this would hold at all integers $x$. In particular, this proposition contains the classical and well-known fact that $\exp(x)$ is not algebraic.

\medskip 
Now we are ready to prove Theorem~\ref{not_dfinite}: we will show that $h_q(x)$ is D-finite (and hence algebraic) if and only if $q$ is a root of unity. This answers Aissen's question completely.

\begin{proof}[Proof of Theorem~\ref{not_dfinite}]
  We already observed that if $q$ is a root of unity, the series $h_q(x)$ is algebraic and hence D-finite. Therefore, one direction is clear and we assume now that $q\in\CC \setminus\{ 0 \}$ is not a root of unity. 
  
  Assume by contradiction that $h_q(x) = \sum_{j\geq0}u_j x^j$ is D-finite. Then $(u_j)_{j\geq0}$ is P-recursive and there exist a positive integer $r$ and $c_0(x),\dots,c_r(x) \in \CC[x]$ with $c_0(x)c_r(x) \neq 0$ such that 
    \begin{equation}\label{rec2}
        u_{j+r}c_r(j) + \dots + u_{j}c_0(j) = 0, \quad \text{for all} \; j\geq0. 
    \end{equation}
    A simple computation shows that 
    \[
        u_{j+1} = u_{j}\frac{\prod_{\ell=1}^a (q^{n+aj+\ell}-1)}{\prod_{\ell=1}^b (q^{k+bj+\ell}-1)\prod_{\ell=1}^{a-b} (q^{n-k+(a-b)j+\ell}-1)}.
    \]
    Then it follows by iteration 
    \begin{align*}
        u_{j+i} &= u_{j} \prod_{m=0}^{i-1} \frac{\prod_{\ell=1}^a (q^{n+a(j+m)+\ell}-1)}{\prod_{\ell=1}^b (q^{k+b(j+m)+\ell}-1)\prod_{\ell=1}^{a-b} (q^{n-k+(a-b)(j+m)+\ell}-1)} \\
        &= u_j \frac{\prod_{\ell=1}^{ia} (q^{n+aj+\ell}-1)}{\prod_{\ell=1}^{ib} (q^{k+bj+\ell}-1)\prod_{\ell=1}^{i(a-b)} (q^{n-k+(a-b)j+\ell}-1)}.
    \end{align*}
    Using this, we may rewrite equation~\eqref{rec2} and obtain
    \begin{equation} \label{sum}
        u_j\left( \sum_{i=0}^r c_i(j) \frac{\prod_{\ell=1}^{ia} (q^{n+aj+\ell}-1)}{\prod_{\ell=1}^{ib} (q^{k+bj+\ell}-1)\prod_{\ell=1}^{i(a-b)} (q^{n-k+(a-b)j+\ell}-1)} \right) = 0,
    \end{equation}
    for all integers $j \geq 0$. Note that $u_j\neq 0$, since $q$ is not a root of unity, hence already the sum above is identical to $0$ for all $j \in \mathbb{N}$. We define
    \[
        P_i(y) \coloneqq \prod_{\ell=1}^{ia} (y^{a}q^{n+\ell}-1)\prod_{\ell=ib+1}^{rb} (y^bq^{k+\ell}-1)\prod_{\ell=i(a-b)+1}^{r(a-b)} (y^{a-b}q^{n-k+\ell}-1) \in \CC[y],
    \]
   so that after multiplication with the common denominator, equation~\eqref{sum} implies that $\sum_{i=1}^r c_i(j)P_i(q^j) = 0$. By Proposition~\ref{prop:jqj} we now obtain that $p(x,y) \coloneqq \sum_{i=1}^r c_i(x)P_i(y)$ must be identically $0$. 
   
   We will show however that $p(x,y)$ cannot be the zero polynomial if the integers $n,k,a,b$ are admissible, more precisely if $a>b>0$. Set first $d \coloneqq \max(\deg(c_i(x)),i=0,\dots,r)$ and write $c_i(x) = \sum_{k=0}^d c_{i,k}x^k$ for some $c_{i,k} \in \CC$. Moreover, let $m$ be an integer such that $c_{r,m} \neq 0$ and denote by $p_m(y)$ the coefficient of $x^m$ in $p(x,y)$. We claim that $p_m(y) \neq 0$. We namely have: 
   \[
        p_m(y) = \sum_{i=0}^r c_{i,m}\prod_{\ell=1}^{ia} (y^{a}q^{n+\ell}-1)\prod_{\ell=ib+1}^{rb} (y^bq^{k+\ell}-1)\prod_{\ell=i(a-b)+1}^{r(a-b)} (y^{a-b}q^{n-k+\ell}-1),
   \]
   and the exponent of $y$ of the leading monomial of each summand is $ia^2+(r-i)b^2+(r-i)(a-b)^2$.  Since we assume that $a>b>0$, it follows that this expression is maximal only for $i=r$, and hence the leading monomial of $p_m(y)$ is $c_{r,m}y^{a^2r}q^{ran + ra(ra+1)/2} \neq 0$.
\end{proof}

\smallskip \noindent {\bf Acknowledgements.} We thank the anonymous referees for helpful suggestions. Both authors were supported by {\href{https://specfun.inria.fr/chyzak/DeRerumNatura/}{DeRerumNatura}} ANR-19-CE40-0018. The second author was also supported by the DOC Fellowship (26101) of the \href{https://www.oeaw.ac.at/}{Austrian Academy of Sciences ÖAW}, and  the \href{https://www.fwf.ac.at/en/}{Austrian Science Fund} (P-34765).


\begin{thebibliography}{10}

\bibitem{Aissen79}
M.~Aissen.
\newblock Variations on a theme of {P}\'{o}lya.
\newblock In {\em Second {I}nternational {C}onference on {C}ombinatorial
  {M}athematics ({N}ew {Y}ork, 1978)}, volume 319 of {\em Ann. New York Acad.
  Sci.}, pages 1--6. 1979.

\bibitem{BeBo92}
J.-P. B\'{e}zivin and A.~Boutabaa.
\newblock Sur les \'{e}quations fonctionelles {$p$}-adiques aux
  {$q$}-diff\'{e}rences.
\newblock {\em Collect. Math.}, 43(2):125--140, 1992.

\bibitem{BCLSS07}
A.~Bostan, F.~Chyzak, G.~Lecerf, B.~Salvy, and E.~Schost.
\newblock Differential equations for algebraic functions.
\newblock In {\em I{SSAC} 2007}, pages 25--32. ACM, New York, 2007.

\bibitem{CaMaMcMo79}
M.~Capobianco, S.~Maurer, D.~McCarthy, and J.~Molluzzo.
\newblock A collection of open problems.
\newblock {\em Annals of the New York Academy of Sciences}, 319(1):565--592,
  1979.

\bibitem{Denjoy32}
A.~{Denjoy}.
\newblock {Sur les courbes definies par les \'equations diff\'erentielles \`a
  la surface du tore}.
\newblock {\em {J. Math. Pures Appl. (9)}}, 11:333--375, 1932.

\bibitem{Desarmenien82}
J.~D\'{e}sarm\'{e}nien.
\newblock Un analogue des congruences de {K}ummer pour les {$q$}-nombres
  d'{E}uler.
\newblock {\em European J. Combin.}, 3(1):19--28, 1982.

\bibitem{Devaney89}
R.~L. Devaney.
\newblock {\em An introduction to chaotic dynamical systems}.
\newblock Addison-Wesley Studies in Nonlinearity. Addison-Wesley Publishing
  Company, Advanced Book Program, Redwood City, CA, second edition, 1989.

\bibitem{Flajolet87}
P.~Flajolet.
\newblock Analytic models and ambiguity of context-free languages.
\newblock volume~49, pages 283--309. 1987.
\newblock Twelfth international colloquium on automata, languages and
  programming (Nafplion, 1985).

\bibitem{Polya22}
G.~P\'olya.
\newblock Sur les s\'eries enti\`eres, dont la somme est une fonction
  alg\'ebrique.
\newblock {\em Enseignement Math.}, 22:38--47, 1921/1922.

\bibitem{Polya69}
G.~P\'{o}lya.
\newblock On the number of certain lattice polygons.
\newblock {\em J. Combinatorial Theory}, 6:102--105, 1969.

\bibitem{Ramis92}
J.-P. Ramis.
\newblock About the growth of entire functions solutions of linear algebraic
  {$q$}-difference equations.
\newblock {\em Ann. Fac. Sci. Toulouse Math. (6)}, 1(1):53--94, 1992.

\bibitem{RiRo14}
T.~Rivoal and J.~Roques.
\newblock Hadamard products of algebraic functions.
\newblock {\em J. Number Theory}, 145:579--603, 2014.

\bibitem{Sagan92}
B.~E. Sagan.
\newblock Congruence properties of {$q$}-analogs.
\newblock {\em Adv. Math.}, 95(1):127--143, 1992.

\bibitem{ScSi19}
R.~Sch\"{a}fke and M.~Singer.
\newblock Consistent systems of linear differential and difference equations.
\newblock {\em J. Eur. Math. Soc. (JEMS)}, 21(9):2751--2792, 2019.

\bibitem{Stanley80}
R.~P. Stanley.
\newblock Differentiably finite power series.
\newblock {\em European J. Combin.}, 1(2):175--188, 1980.

\end{thebibliography}
\end{document}